\documentclass[12pt]{article}
\usepackage{amssymb}
\usepackage{amsmath}
\usepackage[usenames]{color}
\usepackage{graphicx}
\usepackage{amsthm}

\numberwithin{equation}{section}

\usepackage{amsmath,amssymb,mathtools}
\usepackage{amsmath,amssymb,mathtools,amsthm}
\usepackage{textcomp}
\usepackage{amsfonts}
\usepackage{graphicx}
\usepackage{enumerate}
\usepackage{color}
\usepackage{times}
\usepackage[utf8]{inputenc}
\usepackage[T1]{fontenc}
\usepackage{mathtools}
\usepackage{geometry}
\usepackage{bm}
\newcommand{\ncom}{\newcommand}
\newcommand{\non}{\nonumber}
\newcommand{\noi}{\noindent}
\ncom{\vone}{\vskip 2ex}
\ncom{\vtwo}{\vskip 4ex}

\newtheorem{remark}{Remark}
\ncom{\vsa}{\vspace{.3cm}}
\ncom{\vsb}{\vspace{.4cm}}

\begin{document}
	\title{A Unified Approach to Beta Moments, Combinatorial Identities, and Random Walks}

	\author{Puja Pandey and Palaniappan Vellaisamy\\
		\small Department of Statistics and Applied Probability \\[-0.8ex]
		\small University of California Santa Barbara \\[-0.8ex] 
		\small Santa Barbara, CA, 93106, USA. \\
		\small \text{Email}: pujapandey@ucsb.edu;  pvellais@ucsb.edu}
	\maketitle
	
	\begin{abstract}
		The study of random walks has been increasingly popular across diverse disciplines such as statistics, mathematics, quantum physics, where they are used to model paths consisting of successive random steps in a mathematical space. A fundamental quantity of interest is the probability that a simple symmetric random walk returns to the origin after 2n steps. In this paper, a unified probabilistic approach is developed that connects the return probabilities in arbitrary dimensions with moment representations. Using this framework, probabilistic proofs of several combinatorial identities involving beta and gamma functions are provided, and new combinatorial identities in general dimensions are derived.
	\end{abstract}
	\section{Introduction}
	A substantial body of work has been devoted to deriving and interpreting, as well as to obtaining new generalizations of several interesting combinatorial identities involving binomial coefficients, beta and gamma functions and hypergeometric functions. For example, a simple convolution identity for the central binomial coefficients is
	\begin{equation}\label{eqn400}
		\sum_{k=1}^{n} \binom{2k}{k} \binom{2n-2k}{n-k} = 4^{n}.
	\end{equation}

	Identities of this type, together with several related convolution and alternating convolution formulas involving central binomial coefficients, have been studied extensively in the literature \cite{nagy2012combinatorial, mikic2016proof}. Various approaches based on generating functions \cite{charalambides2018enumerative}, combinatorial arguments \cite{de2006pairings, sved1984counting}, recurrence relations, and telescoping techniques have been developed to establish and generalize such identities \cite{mikic2016proof}.
	\\
	\\
	Probabilistic methods involving gamma and beta distributions, hypergeometric functions, and normal approximations have proven to be powerful tools, providing elegant and intuitive proofs of identities that traditionally require advanced analytic techniques  \cite{pathak2018simple, chang2011generalization, vellaisamy2019probabilistic}. These approaches not only yield alternative proofs of known identities but also lead naturally to new identities \cite{vellaisamy2026generalization}, such as
	
	\begin{equation}\label{eqn401}
		\sum_{k=0}^{2n} (-1)^k \binom{2n}{k} \binom{2k}{k}\frac{1}{2^k}= \frac {\binom{2n}{n}}{4^n}
	\end{equation}	
	and
	\begin{equation}\label{eqn402}
	\frac{1}{\pi}\sum_{k=0}^{\infty} \frac{(\frac{1}{2},k)^{2} \; \Gamma(n+\frac{1}{2})}{\Gamma(n+k+\frac{3}{2})} = \frac{1}{4^{2n}} \binom{2n}{n}^{2},
	\end{equation}	
	where $(a, m)= a(a+1)\cdots (a+m-1)$ denotes the Pochhammer symbol.
	Notably, the right-hand sides of these identities correspond to the probability that a simple symmetric random walk returns to the origin after 2n steps in one  \cite{polya1921aufgabe, cicuta1999returns} and two dimensions \cite{polya1921aufgabe, spitzer1976principles}, respectively. 
	\\
	\\
	Simple symmetric random walks occupy a central position in probability theory and statistical physics, serving as canonical models for diffusion, recurrence, and path counting. The simple symmetric random walk on $\mathbb{Z}^{d}$ was first systematically studied by P\'olya \cite{polya1921aufgabe}, who established recurrence for $d=1,2$ and transience for $d \geq 3$. Classical results further describe the probability of return to the origin after a fixed number of steps and reveal a deep combinatorial structure underlying lattice paths in one and higher dimensions  \cite{feller1971introduction, spitzer1976principles}. Variants and extensions of classical random walk have also been extensively studied; among them, the so-called elephant random walk introduces long-range memory effects, where each step depends on the entire history of the walk, leading to anomalous diffusion and markedly different asymptotic behavior \cite{kursten2016random}. While such models capture rich non-Markovian dynamics, the present work focuses exclusively on the memoryless symmetric random walk, where exact return probabilities admit closed-form expressions and, as shown here, admit a natural interpretation in terms of moments of beta-distributed random variables.
	\\
	\\
In deriving the results \eqref{eqn401} and \eqref{eqn402}, moments of both dependent and independent gamma random variables with scale parameter 1 and shape parameter $1/2$ were examined \cite{vellaisamy2026generalization}. This observation naturally motivates the investigation of a deeper connection between moments of certain random variables and return probabilities of random walks, and whether an analogous relationship can be established in three and higher dimensions.
	\\
    \\
	In the course of this investigation, a noteworthy observation was that the 2n-th moment of the random variable $V_{1}=2Y_{1}-1$, where $Y_{1} \sim Be(\frac{1}{2},\frac{1}{2})$ coincides with the probability of return to the origin after 2n steps for a simple symmetric random walk in one dimension. Furthermore, in two dimensions, the probability of return to the origin after 2n steps can be expressed as the 2n-th moment of $(V_{1}+V_{2})/2$, where $V_{1}=2Y_{1}-1$ and $V_{2}=2Y_{2}-1$, with $Y_{1}$ and $Y_{2}$ independent and distributed as $Be(\frac{1}{2},\frac{1}{2})$. Motivated by this connection, a systematic extension of the result to higher dimensions is considered.
	\\
	\\
	The return probabilities can be expressed via multinomial decompositions of lattice paths, \cite{feller1971introduction} and also via Fourier-analytic representations of the lattice Green function \cite{spitzer1976principles}. Extending this framework further, in three dimensions, the probability of return to the origin after 2n steps is given by the 2n-th moment of $(V_{1}+V_{2}+V_{3})/3$, where $V_{i}=2Y_{i}-1, 1 \leq i \leq 3$ and $Y_{i}$'s are  independent and identically distributed as $Be(\frac{1}{2},\frac{1}{2})$ random variables. More generally, in $k$ dimensions the probability of return to the origin after $2n$ steps is given by the $2n$-th moment of $(V_{1}+\cdots+V_{k})/k$, with $V_{i}=2Y_{i}-1$, where $Y_{i} \sim Be\!\left(\frac{1}{2},\frac{1}{2}\right)$ i.i.d.
	\\
	\\
This unified moment-based framework establishes a direct and striking connection between beta-distributed random variables and return probabilities of simple symmetric random walks in one, two, three and higher dimensions. As a byproduct, several new combinatorial identities arising from this probabilistic approach are also obtained, in addition to numerous previously known identities \cite{gould2010combinatorial}.

\section{Beta Moments and Combinatorial Identities}

	Let $X$ be a random variable (r.v.) following symmetric beta distribution $Be(p,p)$, with density
	\begin{equation}\label{eqn100}
		f_{X}(x) = \frac{1}{B(p,p)} x^{p-1} (1-x)^{p-1}, \, \, \, \, \, \, 0 \leq x \leq 1,
	\end{equation}	
	where $B(p,p)$ denotes the usual beta function. It is well known that 
	\begin{equation}\label{eqn600}
		E[X]^{n} = \frac{B(n+p,p)}{B(p,p)}, ~~ n \ge 1.
	\end{equation}	
	Define now $U= 2X-1$. Then, for an integer $ n \ge 1$,
	\begin{align}\label{eqn1}
		E[U]^{2n} &= E[2X-1]^{2n} \nonumber \\
		&= \frac{1}{B(p,p)} \int_{0}^{1}(2x-1)^{2n} x^{p-1} (1-x)^{p-1} du \nonumber \\
	 &= \frac{1}{2 B(p,p)} \int_{-1}^{1} y^{2n} \left(\frac{1+y}{2}\right)^{p-1} \left(\frac{1-y}{2}\right)^{p-1} dy \quad \left( \text{Put} \, 2x-1 = y \right) \nonumber \\
		&= \frac{1}{B(p,p) 2^{2p-2}} \int_{0}^{1} y^{2n} \left(1-y^{2}\right)^{p-1} dy  \nonumber \\
		 &=  \frac{1}{B(p,p) 2^{2p-1}} \int_{t=0}^{1} t^{n-\frac{1}{2}} \left(1-t\right)^{p-1} dt \quad \left( \text{Put} \, y^{2} = t \right) \nonumber \\
		&= \frac{1}{B(p,p) 2^{2p-1}} B(n + \frac{1}{2}, p).
	\end{align}
	 Note that the density of $U = 2X-1$ is
	\begin{align}\label{eqn17}
		f_{U}(u) &= \frac{1}{2^{2p-1} B(p,p)} (1-u^{2})^{p-1}, \, \, \, -1 \leq u \leq 1.
	\end{align}	
Since the density is symmetric about origin, all odd moments vanish, that is, $E[U^{2n+1}] = 0$ for all $n \geq 1$. \\
\\
Let $X_{i}, 1 \leq i \leq k$, be identical and independent random variables following symmetric beta $Be(p,p)$ distributions with density given in \eqref{eqn100}.  \\
Define $U_{i} = 2X_{i} - 1; \, \, 1 \leq i \leq k$, so that the density $U_{i}$ is given in \eqref{eqn17}. We first consider the even moments of a general linear combination $c_1U_1
+c_2U_2 + \ldots + c_kU_k.$
 This computation leads to a unified combinatorial identity that serves as a master formula encompassing several special cases, which have connection to symmetric random walks.
 Let, for simplicity, $c_{(k)}= \sum_{i=1}^{k}c_i.$  Then, using the multinomial expansion, we have
	\begin{align}\label{eqn40}
		 E\left[ \sum_{i=1}^{k} c_{i} U_{i} \right]^{2n} &= E\left[  \sum_{\substack{j_{1}, j_{2},..., j_{k+1} \ge 0 \\ j_{1}+j_{2}+...+j_{k+1} = 2n}} \binom{2n}{j_{1},...,j_{k+1}} c_{(k)}^{j_{1}} \prod_{s=1}^{k} (-2c_{s}X_{s})^{j_{s+1}} \right] \non \\
		 &= \sum_{\substack{j_{1}, j_{2},..., j_{k+1} \ge 0 \\ j_{1}+j_{2}+...+j_{k+1} = 2n}} \binom{2n}{j_{1},...,j_{k+1}} c_{(k)}^{j_{1}} \prod_{s=1}^{k} \Big( (-2c_{s})^{j_{s+1}} E\left[ X_{s}^{j_{s+1}} \right] \Big) \non \\
		 &=  \frac{1}{B(p,p)^{k}}  \sum_{\substack{j_{1}, j_{2},..., j_{k+1} \ge 0 \\ j_{1}+j_{2}+...+j_{k+1} = 2n}} \binom{2n}{j_{1},...,j_{k+1}} c_{(k)}^{j_{1}} \prod_{s=1}^{k} \Bigg( (-2c_{s})^{j_{s+1}} B(j_{s+1}+p,p) \Bigg),
	\end{align}
using \eqref{eqn600}. \\
	\\
Alternatively, one can also compute the moments directly using the moments of $U_{j}$,
	\begin{align}\label{eqn42}
	 E\left[ \sum_{i=1}^{k} c_{i} U_{i} \right]^{2n} &= \sum_{\substack{i_{1}, i_{2},..., i_{k} \ge 0 \\ i_{1}+i_{2}+...+ i_{k} = 2n}} \binom{2n}{i_{1},i_{2},..., i_{k}} \prod_{j=1}^{k}E[c_{j}U_{j}]^{i_{j}}.
	\end{align}	
	Since the odd moments of $U_{i}$ vanish, we have from \eqref{eqn1},
	\begin{align}\label{eqn43}
		E\left[ \sum_{i=1}^{k} c_{i} U_{i}  \right]^{2n} &=   \sum_{\substack{i_{1}, i_{2},..., i_{k} \ge 0 \\ i_{1}+i_{2}+...+ i_{k} = n}} \binom{2n}{2i_{1},2i_{2},...,2i_{k}} \prod_{j=1}^{k}E[c_{j}U_{j}]^{2i_{j}} \non \\
		&=  \sum_{\substack{i_{1}, i_{2},..., i_{k} \ge 0 \\ i_{1}+i_{2}+...+ i_{k} = n}} \binom{2n}{2i_{1},2i_{2},...,2i_{k}} \prod_{j=1}^{k} c_{j}^{2i_{j}} \frac{B(i_{j} + \frac{1}{2},p)}{B(p,p) 2^{2p-1}} \non \\
		&= \frac{1}{B(p,p)^{k} 2^{(2p-1)k}} \sum_{\substack{i_{1}, i_{2},..., i_{k} \ge 0 \\ i_{1}+i_{2}+...+ i_{k} = n}} \binom{2n}{2i_{1},2i_{2},...,2i_{k}} \prod_{j=1}^{k} \Bigg(c_{j}^{2i_{j}} B(i_{j} + \frac{1}{2},p) \Bigg).
	\end{align}
	From \eqref{eqn40} and \eqref{eqn43}, we have a new and general combinatorial identity
	(cancelling out $B^k(p, p)$),
	\begin{align}\label{eqn45}
	 & \sum_{\substack{j_{1}, j_{2},..., j_{k+1} \ge 0 \\ j_{1}+j_{2}+...+j_{k+1} = 2n}} \binom{2n}{j_{1},...,j_{k+1}} c_{(k)}^{j_{1}} \prod_{s=1}^{k} \Bigg( (-2c_{s})^{j_{s+1}} B(j_{s+1}+p,p) \Bigg)    \nonumber  \\
		& \quad \quad \, \, \, \, = \frac{1} {2^{(2p-1)k}} \sum_{\substack{i_{1}, i_{2},..., i_{k} \ge 0 \\ i_{1}+i_{2}+...+ i_{k} = n}} \binom{2n}{2i_{1},2i_{2},...,2i_{k}} \prod_{j=1}^{k} \Bigg(c_{j}^{2i_{j}} B(i_{j} + \frac{1}{2},p) \Bigg) , \, \,  \text{for all} \, \, c_{j} > 0.
	\end{align}	
	\begin{remark}
		If $c_{1}=c_{2}=\cdots=c_{k}=c$ with $c>0$, then $c_{(k)} = kc$ and so \eqref{eqn45} reduces to a symmetric form that is independent of the value of $c$. The resulting expression no longer involves any constant.
			\begin{align}\label{eqn777}
			& \sum_{\substack{j_{1}, j_{2},..., j_{k+1} \ge 0 \\ j_{1}+j_{2}+...+j_{k+1} = 2n}} \binom{2n}{j_{1},...,j_{k+1}} k^{j_{1}} \prod_{s=1}^{k} \Bigg( (-2)^{j_{s+1}} B(j_{s+1}+p,p) \Bigg)    \nonumber  \\
			& \quad \quad \, \, \, \, = \frac{1} {2^{(2p-1)k}} \sum_{\substack{i_{1}, i_{2},..., i_{k} \ge 0 \\ i_{1}+i_{2}+...+ i_{k} = n}} \binom{2n}{2i_{1},2i_{2},...,2i_{k}} \prod_{j=1}^{k} \Bigg( B(i_{j} + \frac{1}{2},p) \Bigg).
		\end{align}	
	\end{remark}

	 In particular, taking $c=1$ yields the same combinatorial identity \eqref{eqn777} independent of any constant. All weight factors disappear, and the combinatorial identity becomes a relation involving only multinomial coefficients and beta functions. 
\\
\\
	In the following subsection, we demonstrate how this general identity simplifies under specific choices of the coefficients, producing well-defined and interpretable combinatorial identities. These identitites when $p = \frac{1}{2}$, they have some connection to random walks.
	\subsection{Some Special Cases}
	We now examine three important special cases of the general identity obtained above. These cases correspond to reductions in dimension and already exhibit rich combinatorial structure. \\
	\\
	\textbf{Case 1:} $C_{s} = 0$, for all $s \geq 4$.  \\
	Setting $c_{s}=0$ for all $s \geq 4$ restricts the problem to a three-variable setting, and only the coefficients $c_{1}, c_{2}, c_{3}$ are retained. The resulting expression captures the essential combinatorial structure of the general formula in a tractable three-dimensional framework. Using $c_{s} = 0$ for all $s \geq 4$, \eqref{eqn45} would reduce to,
\begin{align}\label{eqn567}
& \sum_{\substack{j_{1}, j_{2}, j_{3}, j_{4} \ge 0 \\ j_{1}+j_{2}+j_{3}+j_{4} = 2n}} \binom{2n}{j_{1}, \, j_{2}, \, j_{3}, \, j_{4}} c_{(3)}^{j_{1}}  (-2)^{j_{2}+j_{3}+j_{4}}c_{1}^{j_{2}} c_{2}^{j_{3}} c_{3}^{j_{4}} B(j_{2}+p,p) B(j_{3}+p,p) B(j_{4}+p,p) \nonumber \\
	&= \frac{1}{ (2^{2p-1})^{3}} \sum_{\substack{i_{1}, i_{2},i_{3} \ge 0 \\ i_{1}+i_{2}+i_{3} = n}} \binom{2n}{2i_{1}, \, 2i_{2}, \, 2i_{3}} c_{1}^{2i_{1}} c_{2}^{2i_{2}} c_{3}^{2i_{3}} B(i_{1}+ \frac{1}{2},p)  B(i_{2}+ \frac{1}{2},p) B(i_{3}+ \frac{1}{2},p) .
\end{align}	
This identity represents the complete three-component form of the general moment expression. It captures all interactions among the three variables and provides a natural reference point for the simplified cases that follow. In particular, it lays the groundwork for interpreting the identity in terms of three-dimensional random walk behavior, which will be explored in subsequent sections.
\\
\\
	\textbf{Case 2:} $C_{s} = 0$, for $s \geq 3$.  \\
	Setting $C_{3}= 0$ reduces the problem to a two-variable setting, significantly simplifying the general expression while preserving its essential combinatorial features. 
		\begin{align}\label{eqn206}
		& \sum_{\substack{j_{1},j_{2},j_{3}\ge 0 \\ j_{1}+j_{2}+j_{3} = 2n}} \binom{2n}{j_{1}, j_{2}, j_{3}} (-2)^{j_{2}+j_{3}}  \, c_{(2)}^{j_{1}} \, c_{1}^{j_{2}} \, c_{2}^{j_{3}} \, B(j_{2}+p,p) B(j_{3}+p,p) \nonumber \\ 
		& \qquad \qquad \quad = \frac{1}{(2^{2p-1})^{2}} \sum_{\substack{i_{1}, i_{2} \ge 0 \\ i_{1}+i_{2} = n}} \binom{2n}{2i_{1},2i_{2}} c_{1}^{2i_{1}} c_{2}^{2i_{2}} B(i_{1}+\frac{1}{2},p)  B(i_{2}+\frac{1}{2},p).
	\end{align}	
	The resulting identity can be viewed as a two-dimensional analogue of the general moment formula and will later be linked directly to planar random walks. \\
	\\
	\textbf{Case 3:} $C_{s} = 0$, for $s \geq 2$. \\
	Further specialization to $C_{2} = C_{3} = 0$ yields a one-dimensional setting, where the combinatorial identity admits a particularly simple and transparent form.
	\begin{equation}\label{eqn86}
		\sum_{j=0}^{n} \binom{2n}{j}  \, (-2)^{j} \, B(j+p,p)  = \frac{B(n+\frac{1}{2},p)}{2^{2p-1}}. 
	\end{equation}	
	It is interesting to note that the above identity does not involve $c_{1}$. Despite its simplicity, this identity plays a fundamental role and will serve as the bridge between combinatorial moments and classical one-dimensional random walk results.
	\section{Connection to Random Walks}
	We now establish a direct probabilistic interpretation of the combinatorial identities derived in the previous section. In particular, we show that some special cases of these identities naturally arise in the study of symmetric random walks through return probabilities. \\
	Let $V_{i} = 2(Y_{i} - 1)$, where $Y_{i} \sim Be(\frac{1}{2}, \frac{1}{2})$, $1 \leq i \leq k$, that is,$V_i$'s denote the $U_i$'s, when $p=\frac{1}{2}.$
	\subsection{One Dimensional Random Walk}
	Consider now the rhs of the identity in \eqref{eqn86}. Dividing by $Be(p,p)$ (as we cancelled it earlier), we get 
 \begin{align}\label{eqn444}
 		E[U_{1}^{2n}] = \frac{1}{B(p,p) 2^{2p-1}}  B(n+\frac{1}{2},p).
 \end{align}	
 Using the duplication formula, the following identity follows easily
 \begin{equation}\label{eqn888}
 	\frac{\Gamma{(n+ \frac{1}{2})}}{\Gamma{(\frac{1}{2})}} = \frac{\binom{2n}{n} n!}{4^{n}}.
 \end{equation}	
  Let now $p = \frac{1}{2}$. Then we have
   \begin{align}\label{eqn87}
   E[V_{1}^{2n}] &= \frac{B(n + \frac{1}{2}, \frac{1}{2})}{B(\frac{1}{2}, \frac{1}{2})} \nonumber \\
   &= \frac{\Gamma(n +\frac{1}{2})}{\Gamma(\frac{1}{2})} \frac{1}{n!} \nonumber \\
  	&= \frac{\binom{2n}{n}}{4^{n}},
  \end{align}	
  using \eqref{eqn888}.
 \\
 \\
 	Let $\mathbb{Z}= \{0, 1, \ldots,  \}$ denote the set of integers. Consider the simple symmetric random walk on $\mathbb{Z}$, which moves one step to the right or one step to the left, each with probability $1/2$ each. Let $P_{00}^n(1)$ denote the probability that the random walk, starting at the origin, is again at the origin after $n$ steps. Since the walk changes parity with each step, $P_{00}^n(1) =0 \, \, \, \, \text{when n is  odd.}$
 Thus we focus on  $P_{00}^{2n}(1)$, for $ n \ge 1.$
 \\
 \\
 Let  \( N_{2n} \) denote  the number of paths of { length} \( 2n \) that start and end at the origin. For such a path, the number of steps to the right must equal the number of steps to the left. Therefore, any return to the origin in steps must consist of exactly n right steps and n left steps, in some order. The number of such trajectories is simply the binomial coefficient $\binom{2n}{n}$, since we choose which n of the 2n positions are right steps (and the rest are left steps). Since each trajectory has probability $1/2^{2n},$
 we have
 \begin{align}\label{eqn19}
 	P_{00}^{2n}(1) = \binom{2n}{n}\left(\frac{1}{2}\right)^{2n} = \frac{1}{4^{n}} \binom{2n}{n}. 
 \end{align}
 This is precisely the right-hand side of the one-dimensional version of the identity in \eqref{eqn87}.
 Thus, we have,
 \begin{equation}\label{eqn20}
 	E[V_{1}^{2n}] = P_{00}^{2n}(1) = \frac{1}{4^{n}} \binom{2n}{n}.
 \end{equation}	
 It establishes an exact correspondence between the combinatorial moment identity and the probability of return to the origin in one dimension.
 \begin{remark}
 	At $p= \frac{1}{2}$, \eqref{eqn86} reduces to a combinatorial identity
 	 \begin{align}
 		\sum_{j=0}^{n} \, (-1)^{j}  \frac{1}{2^{j}}  \binom{2n}{j} \binom{2j}{j} = \frac{\binom{2n}{n}}{4^{n}}.
 	\end{align}	
 \end{remark}	
\subsection{Two Dimensional Random Walk}
Consider now the rhs of the identity in \eqref{eqn206}. Dividing $B(p,p)^{2}$ (cancelled earlier), and taking $c_{1}=c_{2} = \frac{1}{2}$, we obtain 
		\begin{align}\label{eqn92}
			E\left[ \frac{U_{1} + U_{2}}{2} \right]^{2n} = \frac{1}{B(p,p)^{2} 2^{4p-2}} \sum_{\substack{i_{1}, i_{2} \ge 0 \\ i_{1}+i_{2} = n}} \binom{2n}{2i_{1},2i_{2}} \left( \frac{1}{2} \right)^{2i_{1} + 2i_{2}} B(i_{1}+\frac{1}{2},p)  B(i_{2}+\frac{1}{2},p). 
		\end{align}
		Let $p = \frac{1}{2}$. Then, we have, from \eqref{eqn888}
	\begin{align}\label{eqn890}
				E\left[ \frac{V_{1} + V_{2}}{2} \right]^{2n} &= \frac{1}{B(\frac{1}{2},\frac{1}{2})^{2} \, 2^{2n}}  \sum_{\substack{i_{1}, i_{2} \ge 0 \\ i_{1}+i_{2} = n}} \frac{2n}{(2i_{1})! \, (2i_{2})!} B(i_{1}+\frac{1}{2},\frac{1}{2})  B(i_{2}+\frac{1}{2},\frac{1}{2}) \nonumber \\
			&= \frac{1}{\pi^{2} \, 2^{2n}}  \sum_{\substack{i_{1}, i_{2} \ge 0 \\ i_{1}+i_{2} = n}} \frac{2n}{(2i_{1})! \, (2i_{2})!} \frac{\Gamma(i_{1}+\frac{1}{2}) \Gamma(\frac{1}{2})}{\Gamma(i+1)} \frac{\Gamma(i_{2}+\frac{1}{2}) \Gamma(\frac{1}{2})}{\Gamma(j+1)} \nonumber \\
			&=  \frac{1}{\pi \, 2^{2n}}  \sum_{\substack{i_{1}, i_{2} \ge 0 \\ i_{1}+i_{2} = n}} \frac{2n}{(2i_{1})! \, (2i_{2})!} \frac{\Gamma(i_{1}+\frac{1}{2})}{i_{1}!} \frac{\Gamma(i_{2}+\frac{1}{2})}{i_{2}!} \nonumber \\
			&= \frac{1}{2^{2n}} \sum_{\substack{i_{1}, i_{2} \ge 0 \\ i_{1}+i_{2} = n}} \frac{2n}{(2i_{1})! \, (2i_{2})!} \frac{\binom{2i_{1}}{i_{1}} }{ 4^{i_{1}}} \frac{\binom{2i_{2}}{i_{2}} }{4^{i_{2}}} \nonumber \\
			&= \frac{1}{ 2^{2n}}  \sum_{\substack{i_{1}, i_{2} \ge 0 \\ i_{1}+i_{2} = n}}  \frac{1}{4^{i_{1}+i_{2}}} \binom{2n}{n} \binom{n}{i_{1},i_{2}}^{2} \nonumber \\
			&= \frac{1}{ 4^{2n}}   \binom{2n}{n} \sum_{\substack{i_{1}, i_{2} \ge 0 \\ i_{1}+i_{2} = n}} \binom{n}{i_{1},i_{2}}^{2}	\nonumber \\
			&= \frac{1}{ 4^{2n}}   \binom{2n}{n}^{2}.
		\end{align}	
		Let, henceforth, $\mathbb{Z}= \{0, 1, \ldots,  \}$ denote the set of non-negative integers. Consider the simple symmetric random walk on the integer points in the plane $\mathbb{Z}^2$ which moves left, right, up, or down, with
		probability $1/4$ each. This  process in indeed called the drunkard’s walk. Let $P_{00}^n(2)$ denote the probability that the random walk returns to the origin in $n$ steps. Then,  $P_{00}^n(2)=0$ when $n$ is  odd, as it can not return to the origin in odd number of steps. \\
		\\
		Let  \( N_{2n} \) denote  the number of paths of { length} \( 2n \) that start and end at the origin. For such a path, the number of steps up must be equal to the number of steps down, and the number of steps to the right must be equal to the number of steps to the left. 
		The total number of steps to the left (also to the right) can be any integer \( k, 0 \le k \le n \); in this case, the trajectory must have \( n - k \) steps up and \( n - k \) steps down. 
		So, if the number of steps to the left is \( k \), the total number of trajectories starting and ending at the origin is the multinational coefficient
		\[
		\binom{2n}{k, k, n-k, n-k}.
		\] 
		
		\noi This implies that the total number of paths of length $2n$ is
		\begin{align*}
			N_{2n} =& \sum_{k=0}^{n} \binom{2n}{k, k, n-k, n-k}   \\
			=& \sum_{k=0}^{n} \frac{(2n)!}{(k!)^2 ((n - k)!)^2}   \\
			= & \sum_{k=0}^{n}  \binom{2n}{n} \binom{n}{k} \binom{n}{n - k}. 
		\end{align*}
		Since each of the above trajectory has the probability 
		$ 1/4^{2n},$ we get
		\begin{align}\label{e33}
			P_{00}^{2 n}(2)=&\sum_{k=0}^{n}  \binom{2n}{n} \binom{n}{k} \binom{n}{n - k} \Big(\frac{1}{4} \Big)^{2n}  \\	
			=& \frac{1}{4^{2n}} \binom{2n}{n} \sum_{k=0}^{n} \binom{n}{k} \binom{n}{n - k} \non \\
			=& \frac{1}{4^{2n}} \binom{2n}{n}^2, \non
		\end{align}
		using Vandermonde identity. \\
		Thus, 
		\begin{equation}
			E\left[\frac{V_{1}+V_{2}}{2}\right]^{2n}= P_{00}^{2n}(2) = \frac{1}{4^{2n}} \binom{2n}{n}^2, \, \, \text{for all} \, \, n \ge 1.
		\end{equation}	
		The resulting return probability coincides precisely with the combinatorial identity derived earlier, confirming the consistency of the framework in two dimensions.
		\begin{remark} At $p=2$ and $c_{1} = c_{2} = \frac{1}{2}$, the general identity \eqref{eqn206} simplifies to an interesting new combinatorial identity.
			\begin{equation}\label{eqn208}
			\sum_{\substack{j_{1}, j_{2}, j_{3} \ge 0 \\ j_{1}+ j_{2}+ j_{3} = 2n}} (-1)^{j_{2}+j_{3}} \, \frac{1}{4^{j_{2}+j_{3}}}   \binom{2n}{j_{1},j_{2},j_{3}}  \binom{2j_{2}}{j_{2}} \, \binom{2j_{3}}{j_{3}} = \frac{1}{ 4^{2n}}   \binom{2n}{n}^{2}
		\end{equation}	
		\end{remark}	
	\subsection{Three Dimensional Random Walk}
	Consider the identity in \eqref{eqn567}. Dividing by $B(p,p)^{3}$ (cancelled earlier) the rhs of \eqref{eqn567} and taking $c_{1} = c_{2} = c_{3} = \frac{1}{3}$, we get
	\begin{align}\label{eqn46}
		& E\left[ \frac{U_{1} + U_{2} + U_{3}}{3}\right]^{2n}  \nonumber \\
		&= \frac{1}{B(p,p)^{3} 2^{6p-3}} \sum_{\substack{i_{1}, i_{2},i_{3} \ge 0 \\ i_{1}+i_{2}+i_{3} = n}} \binom{2n}{2i_{1}, 2i_{2}, 2i_{3}} \left(\frac{1}{3}\right)^{2i_{1}+2i_{2}+2i_{3}} B(i_{1}+ \frac{1}{2},p) B(i_{2}+ \frac{1}{2},p)B(i_{3}+ \frac{1}{2},p) 
	\end{align}	
	Let $p = \frac{1}{2}$. Then, we have, from \eqref{eqn888}
	\begin{align}\label{eqn48}
		& E\left[ \frac{V_{1} + V_{2} + V_{3}}{3}\right]^{2n} \non \\
		&=  \frac{1}{3^{2n} \pi^{3}} \sum_{\substack{i_{1}, i_{2},i_{3} \ge 0 \\ i_{1}+i_{2}+i_{3} = n}} \binom{2n}{2i_{1}, 2i_{2}, 2i_{3}} B(i_{1}+ \frac{1}{2},\frac{1}{2}) B(i_{2}+ \frac{1}{2},\frac{1}{2})B(i_{3}+ \frac{1}{2},\frac{1}{2}) \non \\
		&=  \left(\frac{1}{6}\right)^{2n} \sum_{\substack{i_{1}, i_{2},i_{3} \ge 0 \\ i_{1}+i_{2}+i_{3} = n}} \binom{2n}{2i_{1}, 2i_{2}, 2i_{3}}  \, \binom{2i_{1}}{i_{1}} \binom{2i_{2}}{i_{2}} \binom{2i_{3}}{i_{3}} \non \\
		&= \left(\frac{1}{6}\right)^{2n} \sum_{\substack{i_{1}, i_{2},i_{3} \ge 0 \\ i_{1}+i_{2}+i_{3} = n}} \, \frac{(2n)!}{i_{1}!^{2} i_{2}!^{2} i_{3}!^{2}} 
	\end{align}
	 Consider the simple symmetric random walk on the integer lattice $\mathbb{Z}^3$, where at each step the walk moves one unit in one of the six directions: left, right, up, down, forward, or backward, each with probability $\frac{1}{6}$. Let $P_{00}^n(3)$ denote the probability that the random walk, starting at the origin, returns to the origin after $n$ steps. $P_{00}^n(3) = 0$ when $n$ is odd, since returning to the origin requires equal numbers of positive and negative steps along each coordinate axis. Thus, we focus on $P_{00}^{2n}(3)$, for $n \geq 1$.
	\\
	\\
	Let $N_{2n}$ denote the number of paths of length $2n$ that start and end at the origin. For such a path, the number of steps in each positive direction must equal the number of steps in the corresponding negative direction. Let $i_1, i_2, i_3 \geq 0$ denote the number of steps in the positive directions along the three coordinate axes. Then the number of steps in the corresponding negative directions must also be $i_1, i_2, i_3$, respectively, and we must have $i_1 + i_2 + i_3 = n$.
Thus, for a fixed $(i_1, i_2, i_3)$, the total number of such trajectories is given by the multinomial coefficient
	\[
	\binom{2n}{i_1, i_1, i_2, i_2, i_3, i_3} = \frac{(2n)!}{i_{1}!^{2} \, i_{2}!^{2} \, i_{3}!^{2}}.
	\]
and so  the total number of paths of is
	\begin{align*}
		N_{2n} = \sum_{\substack{i_1, i_2, i_3 \geq 0 \\ i_1 + i_2 + i_3 = n}} 
		\frac{(2n)!}{i_{1}!^{2} \, i_{2}!^{2} \, i_{3}!^{2}}.
	\end{align*}
	
	Since each trajectory has probability $\left(\frac{1}{6}\right)^{2n}$, we obtain
	\begin{align}\label{eqn3d}
		P_{00}^{2n}(3) = \left(\frac{1}{6}\right)^{2n}
	\sum_{\substack{i_1, i_2, i_3 \geq 0 \\ i_1 + i_2 + i_3 = n}} 
		\frac{(2n)!}{i_{1}!^{2} \, i_{2}!^{2} \, i_{3}!^{2}}.
	\end{align}
	The resulting return probability is consistent with the corresponding combinatorial identity in three dimensions, extending the agreement observed in the one and two dimensional cases.
	Thus,
	\begin{equation}
		E\left[\frac{V_1 + V_2 + V_3}{3}\right]^{2n}
		= P_{00}^{2n}(3)
		=\left(\frac{1}{6}\right)^{2n} \sum_{\substack{i_{1}, i_{2},i_{3} \ge 0 \\ i_{1}+i_{2}+i_{3} = n}} \, \frac{(2n)!}{i_{1}!^{2} i_{2}!^{2} i_{3}!^{2}} , \quad \text{for all } n \geq 1.
	\end{equation}
	It is quite striking to observe that the three-dimensional version of this identity aligns beautifully with the combinatorial structure of random walks.
	\begin{remark} At $p= \frac{1}{2}$ and $c_{1} = c_{2} = c_{3} = \frac{1}{3}$, \eqref{eqn567} reduces to a new interesting combinatorial identity
			\begin{equation}\label{eqn52}
			\sum_{\substack{j_{1}, j_{2}, j_{3}, j_{4} \ge 0 \\ j_{1}+j_{2}+j_{3}+j_{4} = 2n}} \binom{2n}{j_{1}, j_{2}, j_{3}, j_{4}} \left(-\frac{1}{6}\right)^{j_{2}+j_{3}+j_{4}} \, \binom{2j_{2}}{j_{2}} \binom{2j_{3}}{j_{3}} \binom{2j_{4}}{j_{4}} =  \left(\frac{1}{6}\right)^{2n} \sum_{\substack{i_{1}, i_{2},i_{3} \ge 0 \\ i_{1}+i_{2}+i_{3} = n}} \, \frac{(2n)!}{i_{1}!^{2} i_{2}!^{2} i_{3}!^{2}}. 
		\end{equation}	
	\end{remark}	
		\subsection{Generalization to $k$ dimension}
		Finally, to highlight the connection of \eqref{eqn45} with a k-dimensional symmetric random walk, divide it by $B(p,p)^{k}$ and take $c_{j} = \frac{1}{k}, 1 \leq j \leq k$
		\begin{align}\label{eqn1000}
			&E\left[\frac{U_1 + \cdots + U_k}{k}\right]^{2n} \non \\
			&=  \frac{1}{B(p,p)^{k} 2^{(2p-1)k}}\sum_{\substack{i_{1}, i_{2},..., i_{k} \ge 0 \\ i_{1}+i_{2}+...+ i_{k} = n}} \binom{2n}{2i_{1},2i_{2},...,2i_{k}} \prod_{j=1}^{k} \Bigg(\left(\frac{1}{k}\right)^{2i_{j}} B(i_{j} + \frac{1}{2},p) \Bigg).
		\end{align}	
		Let $p = \frac{1}{2}$. Then, using \eqref{eqn888}, we have, 
			\begin{align}\label{eqn900}
			 E\left[\frac{V_1 + \cdots + V_k}{k}\right]^{2n} 
			&= \frac{1}{\pi^{k} \, k^{2n}} \sum_{\substack{i_{1}, i_{2},..., i_{k} \ge 0 \\ i_{1}+i_{2}+...+ i_{k} = n}} \binom{2n}{2i_{1},2i_{2},...,2i_{k}} \prod_{j=1}^{k} \Bigg( B(i_{j} + \frac{1}{2},\frac{1}{2}) \Bigg) \nonumber \\
			&= \frac{1}{(2k)^{2n}} \sum_{\substack{i_{1}, i_{2},..., i_{k} \ge 0 \\ i_{1}+i_{2}+...+ i_{k} = n}} \binom{2n}{2i_{1},2i_{2},...,2i_{k}} \prod_{j=1}^{k} \binom{2i_{j}}{i_{j}} \nonumber \\
			&= \frac{1}{(2k)^{2n}} \sum_{\substack{i_{1}, i_{2},..., i_{k} \ge 0 \\ i_{1}+i_{2}+...+ i_{k} = n}} \frac{(2n)!}{i_{1}!^{2} i_{2}!^{2}...i_{k}!^{2}}.
		\end{align}		
		Consider the simple symmetric random walk on $\mathbb{Z}^k$, where at each step the walk moves one unit along any of the $k$ coordinate axes in either direction, each with probability $\frac{1}{2k}$. Let $P_{00}^n(k)$ denote the probability that the walk, starting at the origin, returns to the origin after $n$ steps. As before, $P_{00}^n(k) = 0$ when $n$ is odd, so we restrict attention to $P_{00}^{2n}(k)$.
		\\
		\\
		Let $N_{2n}$ denote the number of paths of length $2n$ that start and end at the origin. For such paths, the number of positive and negative steps along each coordinate must be equal. Let $i_1, \ldots, i_k \geq 0$ with
		$i_1 + \cdots + i_k = n.$
		Then, the number of such trajectories is
		\[
		\binom{2n}{i_1, i_1, \ldots, i_k, i_k} = \frac{(2n)!}{i_{1}!^{2} i_{2}!^{2}...i_{k}!^{2}}.
		\]
		Thus,
		\begin{align}
			P_{00}^{2n}(k) = \left(\frac{1}{2k}\right)^{2n}
			\sum_{\substack{i_{1}, i_{2},..., i_{k} \ge 0 \\ i_{1}+i_{2}+...+ i_{k} = n}}
		\frac{(2n)!}{i_{1}!^{2} i_{2}!^{2}...i_{k}!^{2}}.
		\end{align}
		
		Hence,
		\begin{equation}
			E\left[\frac{V_1 + \cdots + V_k}{k}\right]^{2n}
			= P_{00}^{2n}(k) = \left(\frac{1}{2k}\right)^{2n}
		\sum_{\substack{i_{1}, i_{2},..., i_{k} \ge 0 \\ i_{1}+i_{2}+...+ i_{k} = n}}
			\frac{(2n)!}{i_{1}!^{2} i_{2}!^{2}...i_{k}!^{2}}, \, \, \text{for all} \, \, k \geq 1.
		\end{equation}
			Once again, the probability of returning to the origin after $2n$ steps matches the associated moment computation, completing the dimensional hierarchy.
			\begin{remark} Setting $p = \frac{1}{2}$ in \eqref{eqn45} yields another interesting combinatorial identity
					\begin{align}\label{eqn901}
					\sum_{\substack{j_{1}, j_{2},..., j_{k+1} \ge 0 \\ j_{1}+j_{2}+...+j_{k+1} = 2n}} \binom{2n}{j_{1},...,j_{k+1}} \prod_{s=1}^{k} \left[ \left(-\frac{2}{k}\right)^{j_{s+1}} \binom{2j_{s+1}}{j_{s+1}} \right] = \frac{1}{k^{2n}} \sum_{\substack{i_{1}, i_{2},..., i_{k} \ge 0 \\ i_{1}+i_{2}+...+ i_{k} = n}} \frac{(2n)!}{i_{1}!^{2} i_{2}!^{2}...i_{k}!^{2}}.
				\end{align}	
			\end{remark}	
	\section{Conclusion}
In this paper, a unified connection between even moments of linear combinations of symmetric beta-distributed random variables and return probabilities of symmetric random walks is established. A general combinatorial identity for $E[\sum_{i=1}^{k}c_{i}U_{i}]^{2n}$ is derived, from which several new identities follow. 
\\
\\
It is shown that these identities correspond exactly to the return probabilities of symmetric random walks in $k$ dimensions under appropriate choices of the coefficients and probability $\frac{1}{2}$. This reveals a deep and elegant correspondence between combinatorial moment expansions and probabilistic path counting. Beyond establishing the connection between symmetric random walks and beta moments, the framework developed here also yields several new combinatorial identities.

\end{document}